\documentclass{amsart}

\usepackage[english]{babel}
\usepackage{times}
\usepackage[nothing]{algorithm}
\usepackage{algorithmic}

\newcommand{\boldset}{\mathbb}

\newcommand{\BZ}{\boldset Z}

\newcommand{\BR}{\boldset R}

\DeclareMathOperator{\vol}{vol}

\renewcommand{\vec}[1]{#1}

\theoremstyle{definition}
\newtheorem{definition}{Definition}[section]

\newtheorem{theorem}[definition]{Theorem}
\newtheorem{corollary}[definition]{Corollary}
\newtheorem{proposition}[definition]{Proposition}

\newenvironment{algo}[1]
{
\refstepcounter{definition}
\begin{samepage}
\textbf{Algorithm \thesection.\arabic{definition}} #1\hfill
\begin{center}}
{\end{center}\end{samepage}}

\title{Incremental Algorithms for Lattice Problems}

\author{Boris Hemkemeier}
\address{Universit\"at Dortmund, Fakult\"at f\"ur Mathematik, 44221 Dortmund, Germany}
\email{bhemkemeier@gmail.com}

\author{Frank Vallentin}
\address{Centrum voor Wiskunde en Informatica (CWI), Kruislaan 413, 1098 SJ Amsterdam, The Netherlands}
\email{frank.vallentin@gmail.com}
\thanks{The second author was supported by the Netherlands
Organization for Scientific Research under grant NWO 639.032.203 and by
the Deutsche Forschungsgemeinschaft (DFG) under grant SCHU 1503/4-1.}

\begin{document}

\begin{abstract}
In this short note we give incremental algorithms for the following
lattice problems: finding a basis of a lattice, computing the
successive minima, and determining the orthogonal decomposition. We
prove an upper bound for the number of update steps for every
insertion order. For the determination of the orthogonal decomposition
we efficiently implement an argument due to Kneser.
\end{abstract}

\maketitle

\section{Introduction}
\label{sec:introduction}

Many problems in computational geometry permit a natural computation
by an incremental algorithm. Incremental algorithms process only one
object at a time and insert it into a data structure. Most incremental
algorithms follow an abstract framework: After processing a new object
it is inserted into a data structure. It is first located where the
data structure has to be changed (\textit{localization step}). Then
the data structure has to be updated locally (\textit{update step}) in
order to perform the insertion of a new object.

Here we apply the incremental construction paradigm to the design of
lattice algorithms. Let $\vec{v}_1, \ldots, \vec{v}_m$ be vectors
which span a Euclidean space $E$ and let $L = \BZ \vec{v}_1 + \cdots +
\BZ \vec{v}_m$ be the lattice which is generated by these
vectors. Suppose that we want to compute a property of $L$. First, we
compute the property of the lattice $L_1 = \BZ \vec{v}_1$. Then we
check whether $\vec{v}_2 \in L_1$ (localization step). If $\vec{v}_2
\in L_1$, then nothing has to be done. If $\vec{v}_2 \not\in L_1$,
then we compute the property of the lattice $L_2 = L_1 + \BZ
\vec{v}_2$ (update step), etc. In every update step we compute a
lattice basis for the new lattice $L_i$ which is computationally more
expensive than the localization step. Hence, this algorithmic
framework is attractive if the number of update steps is small.

After fixing notation in Section~\ref{sec:notation} we state an upper
bound for the number of update steps for every insertion order in
Section~\ref{sec:chains}. In Section~\ref{sec:algorithms} we give
algorithms for the following lattice problems: an algorithm which
finds a basis of a lattice given by a set of generators, an algorithm
for the computation of the successive minima of a lattice given by a
complete set of generators (a generating set $S$ of a lattice $L$ is
called \textit{complete} if $S$ contains every vector $v \in L
\setminus \{\vec{0}\}$ with $\|v\| \leq \max_{w \in S} \|w\|$), and an
algorithm for determining the orthogonal decomposition of a lattice
given by a complete set of generators.

These considerations result in a simple meta algorithm with practical
impact. It offers a significant performance benefit compared with
straightforward implementations for classical algorithms. For
experimental results see the technical report
\cite{Hemkemeier-Vallentin-98:Decomposition}. This note is a concise
version of this report where we in particular emphasize the
incremental algorithmic framework.

\section{Notation}
\label{sec:notation}

Let $E$ be a $d$-dimensional Euclidean space. Its inner product is
denoted by $(\cdot, \cdot)$ and the associated norm by $\|\cdot\| =
\sqrt{(\cdot, \cdot)}$. The $d$-dimensional unit ball is denoted
by~$B_d$. A point set $L \subseteq E$ is called a \textit{lattice} if
there exist linearly independent vectors $\vec{b}_1, \ldots, \vec{b}_n
\in E$ such that $L = \BZ \vec{b}_1 + \cdots + \BZ \vec{b}_n$. Then,
$(\vec{b}_1, \ldots, \vec{b}_n)$ is called a \textit{basis} of~$L$ and
$n$ is called the \textit{rank} of~$L$. The \textit{volume} of~$L$ is
given by $\vol L = |\det(\vec{b}_1, \ldots, \vec{b}_n)|$.  A lattice
$L' \subseteq E$ is called a \textit{sublattice} of $L$ if $L'
\subseteq L$. If the rank of $L'$ and $L$ is $d$, then the index of
$L'$ in $L$ is $[L : L'] = \vol L'/\vol L$. The \textit{$k$-th
successive minima} $\lambda_k(L)$ is the minimum value $\lambda$ such
that $\lambda B_d$ contains at least $k$ linearly independent lattice
points of~$L$. We will need the following theorem of Minkowski (see
e.g.\ \cite{Gruber-Lekkerkerker-87:Geometry}).

\begin{theorem} 
\label{th:minkowski}
Let $L \subseteq E$ be a lattice of rank $d$. Then
$$
\frac{2^d}{d!} \vol L \leq \lambda_1(L) \lambda_2(L) \cdots
\lambda_d(L) \vol B_d \leq 2^d \vol L.
$$
\end{theorem}

\section{Chains of Sublattices}
\label{sec:chains}

We want to construct a lattice~$L$, which is generated by the vectors
$\vec{v}_1, \ldots, \vec{v}_m$, incrementally. Update steps are
necessary if $\vec{v}_i \not\in \BZ\vec{v}_1 + \cdots +
\BZ\vec{v}_{i-1}$, where $i = 1,\ldots, m$. The next theorem gives an
upper bound for the number of update steps.

\begin{theorem}
\label{th:main}
Let $\vec{v}_1, \ldots, \vec{v}_m \in E$ be vectors which span~$E$ and which
generate the lattice~$L$. Define $B = \max\limits_{i = 1, \ldots, m}
\|\vec{v}_i\|$. Consider the chain of lattices
\begin{equation}
\label{eq:chain1}
\BZ \vec{v}_1 \;\subseteq\; \BZ \vec{v}_1 + \BZ \vec{v}_2 \;\subseteq\; \ldots
\;\subseteq\; \BZ \vec{v}_1 + \BZ \vec{v}_2 + \cdots + \BZ \vec{v}_m.
\end{equation}
Then, in (\ref{eq:chain1}) inequality holds at most $d +
\log_2(d!(B/\lambda_1(L))^d)$ times.
\end{theorem}

\begin{proof}
First we transform (\ref{eq:chain1}) into a new chain of lattices
which are all of full rank~$d$. Choose indices $1 \leq i_1 < i_2 <
\ldots < i_d \leq m$ such that for all $j \in \{1, \ldots, d\}$ the
rank of $\BZ \vec{v}_{i_1} + \cdots + \BZ \vec{v}_{i_j}$ is $j$ and
$i_j$ is minimal with this property. We define the lattice $L' = \BZ
\vec{v}_{i_1} + \cdots + \BZ \vec{v}_{i_d}$ and consider the
transformed chain
\begin{equation}
\label{eq:chain2}
L' = \BZ \vec{v}_{1} + L' \subseteq \BZ \vec{v}_{1} + \BZ \vec{v}_{2} +
L' \subseteq \ldots \subseteq \BZ \vec{v}_1 + \BZ \vec{v}_2 +
\cdots + \BZ \vec{v}_m + L' = L.
\end{equation}
The number of inequalities in (\ref{eq:chain1}) is at most $d$ plus
the number of inequalities in the chain~(\ref{eq:chain2}). Define
$L'_i = \BZ \vec{v}_1 + \cdots + \BZ \vec{v}_i + L'$, where $i = 1,
\ldots, m$. Since we have
\[
\vol L' / \vol L = [L : L'] = \prod_{i = 2}^m [L'_i: L'_{i - 1}],
\]
the number of inequalities in the chain~(\ref{eq:chain2}) is at most
the number of prime factors of $\vol L' / \vol L$ which is at most
$\log_2(\vol L' / \vol L)$. To finish the proof we apply
Theorem~\ref{th:minkowski} to the quotient $\vol L' / \vol L$ and use
the fact $\lambda_d(L') \leq B$.
\end{proof}

An immediate consequence of Theorem~\ref{th:main} is an upper bound
for the size of a minimal generating set.

\begin{corollary}
\label{cor:covering_theorem}
Let~$L \subseteq E$ be a lattice of full rank~$d$.  Let $S \subseteq
L$ be a finite generating set of~$L$. Then there exists a subset
$S'\subseteq S$ which generates $L$ of size at most $d +
\log_2(d!(B/\lambda_1(L))^d)$ where $B = \max_{\vec{v} \in S}
\|\vec{v}\|$.
\end{corollary}

Note that having long vectors in a lattice basis is not avoidable in
general: Conway and Sloane \cite{Conway-95:Lattice} constructed a
lattice in dimension~$11$ which is generated by its $24$~shortest
vectors but in which no set of $11$~shortest vectors forms a basis.

\section{Algorithms}
\label{sec:algorithms}

In this section we propose algorithms for lattice problems which take
advantage of the incremental construction. The first two algorithms
for computing a lattice basis and for computing the successive minima
are straightforward. For the computation of the unique orthogonal
decomposition we develop new ideas based on an argument of Kneser.

\subsection{Lattice Basis}
\label{ssec:latticebasis}

For computing a lattice basis from a \textit{large} set of generators
we use an algorithm for computing a lattice basis from a
\textit{small} set of generators as a subroutine. Such an algorithm is
the LLL algorithm for linearly dependent vectors of Pohst (see e.g.\
\cite{Cohen-93:Course}, Chapter 2.6.4). Buchmann and Pohst
(\cite{Buchmann-Pohst-89:Computing}) showed for (a variant of) this
algorithm that the number of needed arithmetic operations is
$O(d+m)^4\log B)$. For the incremental algorithm the number of
arithmetic operations is linear in $m$.

\medskip

\begin{algo}{Lattice Basis}
\framebox{
\begin{minipage}{12cm}
\begin{algorithmic}
\label{algo:latticebasis}
\smallskip
\REQUIRE Generating system $\vec{v}_1, \ldots, \vec{v}_m \in E$ of the lattice~$L$.
\smallskip
\ENSURE Basis $\vec{b}_1, \ldots, \vec{b}_n$ of~$L$.
\smallskip
 \STATE $n \leftarrow 0$, $L \leftarrow \{\vec{0}\}$.
 \FOR {$i = 1$ to $m$}
   \IF {$\vec{v}_i \not\in L$}
    \STATE Use a subroutine to get~$n$ and a basis
     $\vec{b}_1, \ldots, \vec{b}_n$ of $L + \BZ \vec{v}_i$.
    \STATE $L \leftarrow \BZ\vec{b}_1 + \cdots + \BZ\vec{b}_n$.
   \ENDIF
 \ENDFOR
\end{algorithmic}
\end{minipage}
}
\end{algo}

\subsection{Successive Minima}
\label{ssec:successiveminima}
For computing the successive minima of a lattice our algorithm is
similar to Algorithm~\ref{algo:latticebasis}. However there are a few
important differences: We need a complete generating system (see
Section~\ref{sec:introduction}), the insertion order is no longer
arbitrary, and in every update step it is enough to compute a basis of
a subspace (instead of a lattice). Hence the number of update steps
equals the rank of the lattice.

\medskip

\begin{algo}{Successive Minima}
\framebox{
\begin{minipage}{12cm}
\begin{algorithmic}
\label{algo:successiveminima}
\smallskip
\REQUIRE Complete generating system $S=\{\vec{v} \in L \backslash
\{\vec{0}\} : \|\vec{v}\|\leq B\}$ of~$L$.
\smallskip
\ENSURE Successive minima $\lambda_1(L), \ldots, \lambda_n(L)$ of $L$.
\smallskip
 \STATE Choose $\vec{v} \in S$ with minimal norm, $S \leftarrow
 S\backslash\{\vec{v}\}$.
 \STATE $n \leftarrow 1$, $U \leftarrow \BR\vec{v}$, $\lambda_n(L)
 \leftarrow \|\vec{v}\|$.
 \WHILE {$S \neq \emptyset$}
  \STATE Choose $\vec{v} \in S$ with minimal norm, $S \leftarrow
  S\backslash\{\vec{v}\}$.
  \IF {$\vec{v} \not\in L$}
   \STATE $U \leftarrow U + \BR \vec{v}$.
   \IF {$\dim U > n$}
    \STATE $n \leftarrow n + 1$, $\lambda_n(L) \leftarrow \|\vec{v}\|$.
   \ENDIF
  \ENDIF
 \ENDWHILE
\end{algorithmic}
\end{minipage}
}
\end{algo}

\subsection{Orthogonal decomposition}
\label{ssec:orthogonaldecomposition}

A lattice is called \textit{decomposable} if it can be written as an
orthogonal direct sum of two non trivial sublattices. Eichler
(\cite{Eichler-52:Note}) proved that every lattice can be decomposed
into indecomposable sublattices which are pairwise orthogonal and that
the decomposition is unique up to order of summands.  In
\cite{Kneser-54:Theorie} Kneser gave a constructive and much simpler
proof. In this section we show how one can efficiently implement
Kneser's argument.

We are given a basis $\vec{b}_1, \ldots, \vec{b}_n \in E$ of the
lattice~$L$, a constant~$B$, and a complete generating system $S =
\{\vec{v} \in L \backslash \{\vec{0}\} : \|\vec{v}\| \leq B\}$. We
want to find the number of indecomposable sublattices~$r$, indices
$i_1 = 1 \leq i_2 < \ldots < i_{r} \leq n < n+1 = i_{r+1}$ and a basis
$\vec{b}'_1, \ldots, \vec{b}'_n$ of $L$ such that for every $j \in
\{1,\ldots, r\}$ the vectors $\vec{b}'_{i_{j}}, \ldots,
\vec{b}'_{i_{j+1}-1}$ form a basis of an indecomposable sublattice.

Now we give Kneser's argument.

\begin{definition}
A vector $\vec{v} \in L \backslash\{\vec{0}\}$ is called
\textit{orthogonal decomposable} if there exist $\vec{x}, \vec{y} \in
L \backslash\{\vec{0}\}$ with $\vec{v} = \vec{x} + \vec{y}$ and
$(\vec{x}, \vec{y}) = 0$.
\end{definition}

The orthogonal indecomposable vectors of~$S$ form the verticex set of an
undirected graph~$G = (V, E)$. In~$G$ two vertices $\vec{v}, \vec{w}
\in V$ are adjacent whenever $(\vec{v}, \vec{w}) \neq 0$. We decompose
$V$ into vertex sets $V_1, \ldots, V_r$ of connected components
of~$G$. Then, the orthogonal decomposition of~$L$ is $L = L_1 \perp
\ldots \perp L_r$ where $L_i$ is the lattice generated by $V_i$.

Using standard algorithms from graph theory one can compute the
connected components in time linear in $O(|V| + |E|)$. In the
following we show that in this case it is possible to compute the
connected components in time linear in $O(|V|)$.

O'Meara observed in \cite{OMeara-80:Indecomposable} that for the
procedure above it is not necessary to determine all orthogonal
indecomposable lattice vectors in~$S$. The length decomposable lattice
vectors are enough:

\begin{definition}
A vector $\vec{v} \in L \backslash\{\vec{0}\}$ is called
\textit{length decomposable} if there exist $\vec{x}, \vec{y} \in L
\backslash\{\vec{0}\}$ with $\vec{v} = \vec{x} + \vec{y}$ and
$\|\vec{x}\| \leq \|\vec{v}\|$ and $\|\vec{y}\| \leq \|\vec{v}\|$.
\end{definition}

On basis of this observation we propose the following algorithm. Its
correctness follows from Proposition~\ref{prop:correct}. In what
follows we denote by $\pi_i$ the orthogonal projection of $E$ onto the
subspace spanned by~$L_i$.

\medskip

\begin{algo}{Orthogonal Decomposition of a Lattice}
\framebox{
\begin{minipage}{12cm}
\begin{algorithmic}
\label{algo:orthogonaldecomposition}
\smallskip
\REQUIRE Complete generating system $S=\{\vec{v} \in L \backslash
\{\vec{0}\} : \|\vec{v}\|\leq B\}$ of~$L$.
\smallskip
\ENSURE Indecomposable sublattices $L_i$ with $L = L_1
\perp \ldots \perp L_r$. 
\smallskip
 \STATE Choose $\vec{v} \in S$ with minimal norm, $S \leftarrow
 S\backslash\{\vec{v}\}$.
 \STATE $r \leftarrow 1$, $L_r \leftarrow \BZ \vec{v}$.
 \WHILE {$S \neq \emptyset$}
  \STATE Choose $\vec{v} \in S$ with minimal norm, $S \leftarrow
  S\backslash\{\vec{v}\}$.
  \IF {$\vec{v} \not\in \sum_{i=1}^r L_i$}
   \STATE $J \leftarrow \{ j \in
    \{1, \ldots, r\} : \pi_j(\vec{v}) \neq \vec{0}\}$.
   \STATE $M \leftarrow \BZ \vec{v} + \sum_{i\in J} L_i$.
   \STATE $\{L_1, \ldots, L_{r-|J|}\} \leftarrow \{L_i : i \notin J\}$,
    $L_{r-|J|+1} \leftarrow M$,  $r\leftarrow {r-|J|+1}$.
  \ENDIF
 \ENDWHILE
\end{algorithmic}
\end{minipage}
}
\end{algo}

\medskip

\begin{proposition}
\label{prop:correct}
At the end of each iteration the computed sublattices are indecomposable
and pairwise orthogonal.
\end{proposition}

\begin{proof}
By induction the sublattices $L_1, \ldots, L_r$ are indecomposable and
pairwise orthogonal.  Let $\vec{v}$ be a shortest vector in~$S$.  If
$\vec{v} \not\in \sum_{i=1}^r L_i$, then $\vec{v}$ is not length
decomposable.  In particular we have either $\pi_i(\vec{v}) = \vec{0}$
or $\pi_i(\vec{v}) \notin L_i$ where $i = 1, \ldots, r$. Define $J =
\{j \in \{1, \ldots, r\} : \pi_j(\vec{v}) \neq \vec{0}\}$. One can
choose vectors $\vec{v}_j \in L_j$, where $j \in J$, which are not
length decomposable and which are not orthogonal to~$\vec{v}$. In the
graph $G$ these vectors are all adjacent to $\vec{v}$.  Hence,
$\BZ\vec{v} + \sum_{j \in J} L_j$ is indecomposable and we get
$\sum_{i\in I \backslash J} L_i \perp (\BZ \vec{v} + \sum_{j\in J}
L_j)$ because $\pi_i(\vec{v}) = \vec{0}$ for $i \in I \backslash J$.
\end{proof}

\section{Acknowledgements}

We thank Martin Kneser for pointing out the reference to
\cite{OMeara-80:Indecomposable}.

\end{document}